\documentstyle[11pt,leqno,amscd,amssymb,pstricks,verbatim]{amsart}

\newtheorem{thm}{Theorem}[section]
\newtheorem{lem}[thm]{Lemma}

\newtheorem{prop}[thm]{Proposition}

\theoremstyle{definition}

\newtheorem{defn}[thm]{Definition}

\newtheorem{rem}[thm]{Remark}

\numberwithin{equation}{thm}

\def\cc{{\cal C}}  
   \def\cb{{\cal B}}

\def\ca{{\cal A}}  \def\cu{{\cal U}}
  \def\bbz{{\mathbb Z}}

  \def\leq{\leqslant}  \def\geq{\geqslant}

\def\Hom{\mbox{\rm Hom}}

\def\Ext{\mbox{\rm Ext}\,}   
   \def\End{\mbox{\rm End}\,}
\def\udim{{\mathbf dim}} 
  
  \def\top{\mbox{\rm top}\,}
  \def\soc{\mbox{\rm soc}\,}

  \def\fkg{{\frak g}}

\def\bA{{\mathbf A}}
\begin{document}

\title[The bijection between Exceptional Subcategories And Non-crossing Partitions]%
{The bijection between Exceptional Subcategories And Non-crossing Partitions}

\author{Anningzhe Gao}

\begin{abstract}
This note discusses the bijection between the exceptional subcategories of representations of quivers and generalized non-crossing partitions of Weyl groups.  We give a new proof of the Ingalls-Thomas-Igusa-Schiffler bijection by using the exchange property of the Weyl groups of the Kac-Moody Lie algebras.
\end{abstract}

\maketitle


\section{Introduction}

Representations of quivers have deep relations with the Kac-Moody Lie algebras. Once we are given an acyclic quiver Q, we can define its representation category $repQ$ over a field $k$. It is an abelian category. Let $\bA (mod\Lambda)$ denote the set of all exceptional subcategories of $mod\Lambda$ where $\Lambda=kQ$ the path algebra of Q and $mod\Lambda$ denotes the finite dimensional (left) $\Lambda$-modules. Let $K_0(repQ)$ be the Grothendieck group, the symmetric Euler form $(-,-)$ of Euler form is well-defined in $K_0(repQ)$. The system $\{K_0(repQ), (-,-)\}$ then can define a generalized Cartan matrix. The correspond Kac-Moody Lie algbra is denoted by $\fkg(\Lambda)$. And the weyl group of $\fkg(\Lambda)$ is denoted by $W(\Lambda)$, the Coxeter element is $c(\Lambda)$. We define what we call exceptional subcategories and generalized non-crossing partitions. We consider the following

\textbf{Main Theorem}: There is an isomorphism 
\[cox: \bA(mod\Lambda)\rightarrow Nc(W(\Lambda),c(\Lambda))\] 	

We will define the isomorphism in Section 2. This Theorem was first proved by Ingalls and Thomas \cite{IT} for Dynkin and tame case.Then Igusa and Schiffler \cite{IS} proved it in general case. In this note, we will give an elementary and straightforward proof for this theorem. 

In Section 2 we introduce some basic definitions and preliminary results to give a definition of the map $cox$. In Section 3 we recall the braid group action on exceptional sequences due to Crawley-Boevey\cite{CB1}.
 In Section 4 we show that there is a natural action of braid group on weyl group which is called Hurwitz transform and this action is transitive on the generalized non-crossing partitions. With the preparation of Section 3 and 4, Section 5 gives a proof of the bijection between the exceptional sequences and generalized non-crossing partitions.  


\bigskip

\section{Definitions and preliminary results} An acyclic quiver is an oriented graph Q without oriented cycles. We write it as $Q=(Q_0,Q_1)$, where $Q_0$ is the set of all vertices and $Q_1$ is the set of all arrows. Consider the representation category of  $repQ$ over a field $k$. Let $\Lambda=kQ$ the path algebra of $Q$. There is a canonical category equivalence 
$mod\Lambda\simeq repQ$.  In the paper we will identify
representations of $Q$ over $k$ with $\Lambda$-modules.

A representation of $Q$ is denoted by $(V_i,v_\alpha,i\in Q_0, \alpha\in Q_1)$. Here $v_\alpha$ is a linear transform from $V_{h(\alpha)}$ to $V_{t(\alpha)}$ where $h(\alpha)$ is the head of the arrow $\alpha$ and $t(\alpha)$ is the tail of it.

Let $K_0(mod\Lambda)$ be the Gronthendieck group. Then $K_0(mod\Lambda)$ is a free abelian group of rank n=$\#Q_0$. So $K_0(mod\Lambda)\cong  \bbz^n$. Given $M\in mod\Lambda$ the dimension vector of it is $$\udim(M)=(dimV_1,...,dimV_n)$$
Given two vectors $v,w\in K_0(mod\Lambda)$, the bilinear form Euler form is defined as follows: $$\langle v,w\rangle=\Sigma_{i\in Q_0} v_iw_i-\Sigma_{\alpha\in Q_1}v_{h(\alpha)}w_{t(\alpha)}$$
For two modules $M,N\in mod\Lambda$, we define $\langle M,N \rangle=\langle \udim(M),\udim(N)\rangle$. The symmetric Euler form $(-,-)$  is defined by $(v,w)=\langle v,w\rangle+\langle w,v \rangle$. The system $(K_0(mod\Lambda),(-,-))$ then defines a generalized Cartan matrix and then we obtain the corresponding Kac-Moody Lie algebra $\fkg(\Lambda)$. Let $\Phi$ be its root system. The dimension vectors of simple objects $\{\udim(S_i) \}_{i\in Q_0}$ are exactly the simple roots of $\Phi$. Given an element $v\in K_0(mod\Lambda)$, if we write $v=\Sigma_{i=1}^n c_i\udim S_i$, the support of $v$ is the subset of the bases of $K_0(mod\Lambda)$  such that $c_i\neq 0$. We say $v$ is positive if $c_i>0$ for $\udim S_i$ in its support and $v\neq 0$. Then $\Phi$ has the decomposition $\Phi=\Phi^+\cup\Phi^-$ where $\Phi^+$ is the set of positive roots and $\Phi^-=-\Phi^+$. The real roots is the root that can be obtained from simple roots by reflections. We denote the complement of real roots imaginary roots. So $\Phi=\Phi_{re}\cup\Phi_{im}$. For every real root $v$, the equality $(v,v)=2$ holds.
 With this notion we have the reflection transforms for all real roots as
\[\sigma_{v}(w)=w-(v,w)v\]
For each indecomposable module $M$ such that $\udim(M)$ is a real root, we define $\sigma_M=\sigma_{dimM}$. The equation $\Phi_{re}=\cup_{i\in Q_0}W(\Lambda){\udim(S_i)}$ holds.

For an element $\omega\in W(\Lambda)$, we define its absolute length $\lvert\omega\rvert_a$ equal to the minimal number $l$ that $\omega$ can be written as product of $l$ reflections of real roots. With the absolute length we define a partial order on $W(\Lambda)$ by the following:
$$\omega_1\leq\omega_2\Leftrightarrow\lvert\omega_1\rvert_a+\lvert\omega_1^{-1}\omega_2\rvert_a=\lvert\omega_2\rvert_a$$

If $\{i_1,...,i_n\}=Q_0$, then $\sigma_{S_{i_1}}...\sigma_{S_{i_n}}$ is the Coxeter element in $W(\Lambda)$.
Choose one Coxeter element $c$, we define the set of generalized non-crossing partitions $Nc(W,c)$ as
\[Nc(W,c)=\{\sigma\in W\arrowvert \sigma\leq c\}\]

For $i\in Q_0$, we have simple modules $S_i$, the indecomposable projective modules $P_i$ with $\top P_i=S_i$, the indecomposable injective modules $I_i$ with $\soc I_i=S_i$. Then $\{S_i\}_{i\in Q_0}$ is the complete collection of the simple modules, $\{P_i\}_{i\in Q_0}$ is the complete collection of the indecomposable projective modules, $\{I_i\}_{i\in Q_0}$ is the complete collection of the indecomposable injective modules.

 A module M is called exceptional if $\End_\Lambda(M)=k$ and $\Ext_\Lambda(M,M)=0$. An antichain is a set of modules
 \[\{A_1,A_2,...,A_r\}\]
 such that $\Hom_\Lambda(A_i,A_j)=0$ for all $i\neq j$ and $\Hom_\Lambda(A_i,A_i)=k$. Recall that given an antichain, the Ext-quiver of it is defined as follows: The vertices of the quiver are the elements in the antichain, and there is an arrow from $i$ to $j$ if $\Ext_\Lambda(A_i,A_j)\neq 0$. An antichain is called exceptional if its Ext-quiver is acyclic\cite{DX}.  Then given an exceptional antichain, we can define an exceptional subcategory $\mathcal{A}$ as its extension closure. Then we can show that $\mathcal{A}$ is closed under extension, kernel of monomorphism and cokernel of epimorphism, which we call it a thick subcategory. By $\mathcal{A}$$\leq$$\mathcal{B}$ we mean $\mathcal{A}$$\subseteq$$\mathcal{B}$.
 
 A sequence $E=(E_1,E_2,...E_r)$ in $mod\Lambda$ is called an exceptional sequence if each $E_i$ is an exceptional $\Lambda$ module and we have $\Hom_\Lambda(E_j,E_i)=\Ext_\Lambda(E_j,E_i)=0$ for $j>i$. If $r=n$, we call the exceptional sequence a complete exceptional sequence.
 
 Given an exceptional sequence $E$, we can define a full subcategory $\mathcal{A}$ (denote by $\cc(E)$) of $mod\Lambda$ as the thick closure of the sequence. On the other hand, for every exceptional subcategory $\mathcal{A}$, all its simple object $\{S_1,S_2,...,S_r\}$ is an exceptional antichain, we can relabel it such that $(S_1,...,S_r)$ is an exceptional sequence.

 With the above notions and properties we can define the bijection between the exceptional subcategories and generalized non-crossing partitions.
 
 For every exceptional subcategory $\mathcal{A}$, choose a complete exceptional sequence $E=(E_1,E_2,...,E_r)$, define a correspondence
 $$cox:\ \bA(mod\Lambda)\rightarrow\ Nc(W(\Lambda),c(\Lambda))$$
 $$cox(\ca)=\sigma_{E_{\text{1}}}\sigma_{E_{\text{2}}}...\sigma_{E_{\text{r}}}$$

 In section 3, we will prove that this is a well defined map and in the last section we will prove that this map is actually a bijection.
 
 \section{Braid group action on exceptional sequences}

 Recall that a braid group $B_n$ is a group generated by \{$\rho_1,\rho_2,...,\rho_{n-1}$\} with respect to the following relations:
 
 1) $\rho_i\rho_{i+1}\rho_i=\rho_{i+1}\rho_i\rho_{i+1}$
 
 2) $\rho_i\rho_j=\rho_j\rho_i$ for $\lvert j-i\rvert\geq2$
 
 We introduce some well known lemmas which are taken from \cite{CB1}.
 
 As above,let $Q$ be an acyclic quiver with n vertices. Let $\Lambda=kQ$. First we define the perpendicular subcategory.
 
 \begin{defn}
 	Given a subcategory $\mathcal{U}$ of $mod\Lambda$. The $right$ ($resp.left$) perpendicular subcategory of $\cu$ which is denoted by $\mathcal{U}^{\perp}$ ($resp.^\perp\cu$) the set
 	\[\cu^{\perp}=\{M\in mod\Lambda\arrowvert \Hom_\Lambda(N,M)=\Ext_\Lambda(N,M)=0 \forall  N\in \cu\}\]
 	$$(resp.^{\perp}\cu=\{M\in mod\Lambda\arrowvert \Hom_\Lambda(M,N)=\Ext_\Lambda(M,N)=0 \forall  N\in \cu\})$$
 \end{defn}

 Now we say a pair $(\mathcal{U},\mathcal{V})$ is perpendicular pair if $\mathcal{U}=\mathcal{V}^{\perp}$ and $\mathcal{V}=^{\perp}\mathcal{U}$.
 
 Use these notation, we can describe the following lemmas.

 \begin{lem}
 	If $E=(E_1,E_2,...,E_r)$ is an exceptional sequence, then $\cc(E)^{\perp}$ ( $^{\perp}\cc(E)$) is equivalent to $kQ(E^\perp)-mod$ ($kQ(^\perp E)-mod$) category where $Q(E^\perp)$ ($Q(^\perp E)$) is some acyclic quiver with $(n-r)$ vertices.
 \end{lem} 
 
 \begin{pf}
 	We refer to Schofield's paper [5,Theorem 2.3]. 
 \end{pf}
 \begin{lem}
 	For a complete exceptional sequence $E=(E_1,E_2,..,E_n)$, we have $\cc (E)=mod\Lambda$.
 \end{lem} 
 \begin{pf}
    We prove this lemma by induction on the number of vertices of $Q$. When n=1, it is easy to see. Now we suppose for $k<n$, the lemma holds.
    	
    Let $X=E_n$. Now $E'=(E_1,...,E_{n-1})$ is a complete sequence of $kQ$($X^\perp$). By induction we have $\cc(E')$=$X^\perp$. 
    
    Suppose that $X$ is not a projective module. Then by Bongartz completion we have $Y\in X^\perp$ such that $T=X\oplus Y$ is a tilting module. Since by the definition of a tilting module, there is an exact sequence $$0\rightarrow\varLambda\rightarrow T'\rightarrow T''\rightarrow0$$ where $T', T''\in add(T)$, we can conclude that all projectives are in $\cc(E)$. Since every module has a projective resolution, the lemma has been proved.
    
    If $X$ is projective, let $X=P(i)$ for some i. Then $X^\perp$ is just the category of the representations of the quiver $Q'$ which is obtained by deleting vertex $i$ in $Q$. Since we have the exact sequence $0\rightarrow radX\rightarrow X\rightarrow S_i\rightarrow0$, $S_i\in\cc(E)$. For $j\neq i$, $S_j\in\cc(E')$ by induction. So all the simple modules are in the $\cc(E)$. We finish the proof.
 \end{pf}
 
 \begin{lem}
 	Each exceptional sequence $E=(E_1,E_2,..,E_r)$ can be extended to a complete exceptional sequence. And for  exceptional subcategory $\cu$, we have 
 	\[^\perp(\mathcal{U}^\perp)=(^\perp\mathcal{U})^\perp=\mathcal{U}\]
 	
 \end{lem} 
 \begin{pf}
 	Since $E=(E_1,...,E_r)$ is an exceptional sequence, from Lemma 3.2, $\cc(E)^\perp$ is equivalent to the representation category of an acyclic quiver of $n-r$ vertices. So we can choose a complete exceptional sequence $F$ of $kQ(E)-mod$. Then $(F,E)$ is a complete exceptional sequence of $mod\Lambda$. For $^\perp \cc(E)$, things are similar. The first statement is proved. 
 	
 	For the second statement, it is obvious that $(^\perp\mathcal{U})^\perp\subseteq\mathcal{U}$. We already knew that there is a complete exceptional sequence having the form $(E,F)$ in $mod\Lambda$ where $E$ is a complete sequence of $\mathcal{U}$ and $F$ is the complete sequence of $^\perp \cc(E)$. Then by Lemma 3.3 $$^\perp \cc(E)=\cc(F)$$
 	So $(^\perp\mathcal{U})^\perp$=$\cc(F)^\perp\supseteq \cc(E)$
 	Then the second statement is proved.
 \end{pf}
 
 \begin{lem}
 	If $E=(E_1,E_2,...,E_{i-1},X,E_{i+1},...,E_n)$ 
 	and 
 	
 	$E'=(E_1,E_2,...,E_{i-1},Y,E_{i+1},...,E_n)$ both are exceptional sequences, then $X\cong Y$.
 	
 \end{lem} 
 
\begin{pf}
    By passing to $^\perp(E_1,E_2,...,E_{i-1})$ and 
	
	$(E_{i+1},...,E_n)^\perp$, we obtain an exceptional subcategory with only one simple object . So $X\cong Y$.
\end{pf}

The following lemma is due to Schofield which is well known, for proof, see\cite{SCH}.
 \begin{lem}
 	For any exceptional module $M$, if $M$ is not simple in $mod\Lambda$, then there exists two exceptional modules $X,Y$ such that $\Hom_{\Lambda}(X,Y)=\Hom_{\Lambda}(Y,X)=\Ext_{\Lambda}(Y,X)=0$ and $M$ is relative project in $\cc(X,Y)$ and there exists an exact sequence
 	\[0\rightarrow Y^b\rightarrow M\rightarrow X^a\rightarrow0\].
 \end{lem}

 \begin{lem}
 	For any exceptional pair $(X,Y)$, there exists a unique exceptional module $R_Y X\in \cc(X,Y)$$(resp. L_X Y\in \cc(X,Y))$ such that $(Y,R_Y X)(resp.(L_X Y,X))$ is an exceptional pair.
 \end{lem} 
 \begin{pf}
 	Since $\cc(X,Y)$ can be viewed as the representation category of an acyclic quiver with 2 vertices, by Lemma 3.4, $X$ ($Y$) can be extended from the left (right) to a complete exceptional sequence of $\cc(X,Y)$. Then we finish the proof.
 \end{pf}
 
 Now we can introduce the braid group action on the complete exceptional sequences.
 
 \begin{defn}
 	Given a complete exceptional sequence $E=(E_1,E_2,...,E_n)$, we define the braid group actions as follows:
 	\[\rho_i(E_1,E_2,...,E_n)=(E_1,..,E_{i-1},E_{i+1},R_{E_{i+1}}E_i,E_{i+2},..,E_n)\] 
 	\[\rho_i^{-1}(E_1,E_2,...,E_n)=(E_1,..,E_{i-1},L_{E_i}E_{i+1},E_i,E_{i+2},..,E_n)\]
 \end{defn} 
 
 We can check by calculation directly that this is a $B_n$ action on the complete exceptional sequences. Then Crawley-Boevey proved that this $B_n$ action is transitive:
 
 \begin{thm}
 	The $B_n$ action on the set of complete exceptional sequences is transitive.
 \end{thm} 

 Since the antichain of an exceptional subcategory is a complete exceptional sequence, we have the following proposition, see\cite{Ringel}.

 \begin{prop}
 	For any two complete sequences of an exceptional subcategory $\mathcal{A}$: $E=(E_1,...,E_r)$ and $E'=(E'_1,...,E'_r)$, we have $\sigma_{E_1}\sigma_{E_2}...\sigma_{E_r}=\sigma_{E'_1}\sigma_{E'_2}...\sigma_{E'_r}$
 \end{prop} 
 
 \begin{pf} By Theorem 3.9, the braid group acts transitively on the set of complete sequences. And we have the formulas
 	\[\udim R_YX=\pm\sigma_Y(X)\]
 	\[\udim L_XY=\pm\sigma_X(Y)\]
 	(Crawley-Boevey's paper)\cite{CB1}
 	
 	Fix a complete exceptional sequence $E=(E_1,E_2,...,E_r)$. First we prove that for a generator $\rho_i$ of $B_n$ and denote $E^*=(E^*_1,E^*_2,...,E^*_r)=\rho_i E$, the equation $\sigma_{E_1}\sigma_{E_2}...\sigma_{E_r}=\sigma_{E^*_1}\sigma_{E^*_2}...\sigma_{E^*_r}$ holds. 
 	
 	By definition, $$\sigma_{E^*_1}\sigma_{E^*_2}...\sigma_{E^*_r}=\sigma_{E_1}...\sigma_{E_{i+1}}\sigma_{L_{E_i}E_{i+1}}...\sigma_{E_r}=\sigma_{E_1}...\sigma_{E_{i+1}}\sigma_{\sigma_{E_i}E_{i+1}}...\sigma_{E_r}=\sigma_{E_1}\sigma_{E_2}...\sigma_{E_r}$$
 	
 	So we conclude that $B_n$ action does not change the product.
 	 
 	Then for every two complete sequences $E=(E_1,...,E_r)$ and $E'=(E'_1,...,E'_r)$, we have $\sigma_{E_1}\sigma_{E_2}...\sigma_{E_r}=\sigma_{E'_1}\sigma_{E'_2}...\sigma_{E'_r}$ because braid group action on the set of complete exceptional sequences is transitively.
 	
 	Thus the lemma is proved. 
 \end{pf}

\section{Braid group action on the set of non-crossing partitions}

As we introduced in Section 2, from an acyclic quiver we can get a Kac-Moody Lie algebra. The aim of this section is to prove that Hurwitzs transformation is transitive in $Nc(W,c)$.

By $\Phi$ and $W$ we denote the roots system and Weyl group of the Kac-Moody Lie algebra $\fkg(\Lambda)$ respectively. Let  $\{S_1,S_2,...,S_n\}$ be the complete collection of non-isomorphic simple modules of $mod\Lambda$. There are natural decompositions $\Phi=\Phi^+\cup\Phi^-$ and $\Phi=\Phi_{re}\cup\Phi_{im}$ as we discussed in Section 2. According to Kac's theorem, $\Phi^+$ equal to theset of the dimension vectors of indecomposable modules in $mod\Lambda$. Moreover, $\Phi^-=-\Phi^+$. We write $\alpha>0$ for an element $\alpha$ in $K_0(mod\Lambda)$ if $\alpha\neq 0$ and $\alpha=\Sigma_{i=1}^n c_i[S_i]$ where $c_i\in\bbz_{\geq 0}$ in $K_0(mod\Lambda)$. We write $\alpha>\beta$ if $\alpha-\beta>0$.

It is well known that $W(\Phi_{re})=\Phi_{re}$ and $W(\Phi_{im})=\Phi_{im}$

\begin{lem} The simple reflection $\sigma_{S_i}$ preserves $\Phi^+-\{\udim S_i\}$
\end{lem}

\begin{pf}
	If $\alpha\in\Phi^+$ and $\alpha\neq \udim S_i$, then $\sigma_{S_i}(\alpha)$ can not be a negative root. Since $\sigma_{S_i}$ transforms roots to roots, then $\sigma_{S_i}(\alpha)\in\Phi^+$.
	\end{pf}
	
The following lemma is the well-known exchange property.

\begin{lem} If  $\sigma_{S_{i_1}}\sigma_{S_{i_2}}..\sigma_{S_{i_k}}(\alpha)<0$ for some $\alpha\in\Phi^+\cap\Phi_{re}$, then there is $1\leq t\leq k$ such that $\sigma_{S_{i_t}}\sigma_{S_{i_{t+1}}}..\sigma_{S_{i_k}}=\sigma_{S_{i_{t+1}}}\sigma_{S_{i_{t+2}}}..\sigma_{S_{i_k}}\sigma_{\alpha}$
\end{lem}

\begin{pf} 
	Since $\alpha>0$ there exists $1\leq t\leq k$ such that $\sigma_{S_{i_{t+1}}}\sigma_{S_{i_{t+2}}}..\sigma_{S_{i_k}}(\alpha)>0$ and $\sigma_{S_{i_{t}}}\sigma_{S_{i_{t+1}}}..\sigma_{S_{i_k}}(\alpha)<0$. By Lemma 4.1, the equality $\sigma_{S_{i_{t+1}}}\sigma_{S_{i_{t+2}}}..\sigma_{S_{i_k}}(\alpha)=\udim S_{i_t}$ holds. So
\[\sigma_{S_{i_t}}=\sigma_{S_{i_{t+1}}}\sigma_{S_{i_{t+2}}}..\sigma_{S_{i_k}}\sigma_\alpha(\sigma_{S_{i_{t+1}}}\sigma_{S_{i_{t+2}}}..\sigma_{S_{i_k}})^{-1}\]
The lemma follows.
\end{pf}

We now give the concrete definition of the absolute length (See Section 2).

\begin{defn}
	An element $\omega\in W$ has an absolute length $\lvert \omega\rvert_a$ if $\omega$ can be written as products of $\lvert \omega\rvert_a$ reflections but can not be written by product of less number of reflections.
	
\end{defn}

Let $T$ be the set of all reflections at real roots in $W$. The following defines the braid group action on $T^n$, called Hurwitz transformation. The definition of the braid group $B_n$ is already given at the beginning of Section 3.

\begin{defn}
	Given $(\sigma_{\alpha_1},\sigma_{\alpha_2},...,\sigma_{\alpha_n})\in T^n$, the Hurwitz transformation on $T^n$ is defined by for the canonical generators $\rho_i$ of $B_n$:
	\[\rho_i(\sigma_{\alpha_1},...\sigma_{\alpha_i},\sigma_{\alpha_{i+1}}...,\sigma_{\alpha_n})=(\sigma_{\alpha_1},...,\sigma_{\alpha_{i+1}},\sigma_{\sigma_{\alpha_{i+1}}(\alpha_i)},...,\sigma_{\alpha_n})\]
\end{defn} 

It can be checked directly by calculation that this is a group action of $B_n$ on $T^n$. 

\begin{rem}
	From the definition of the action, what should be noticed is that $$\sigma_{\alpha_{i+1}}\sigma_{\sigma_{\alpha_{i+1}}(\alpha_i)}=\sigma_{\alpha_{i+1}}\sigma_{\alpha_{i+1}}\sigma_{\alpha_{i}}\sigma_{\alpha_{i+1}}=\sigma_{\alpha_{i}}\sigma_{\alpha_{i+1}}$$
	so the  action of $B_n$ on $T^n$ does not change the product of $(\sigma_{\alpha_1},\sigma_{\alpha_2},...,\sigma_{\alpha_n})$
	Thus it induces an action of $B_n$ on $Nc(W,c)$.
\end{rem}

We now label the simple objects of $mod\Lambda$ in an appropriate order such that $S=(S_1,S_2,...,S_n)$ is a complete exceptional sequence.

\begin{thm}
	If $\sigma_{\alpha_1}\sigma_{\alpha_2}...\sigma_{\alpha_n}=\sigma_{S_1}\sigma_{S_2}...\sigma_{S_n}$ where all $\alpha_i$ are positive real roots, then $(\sigma_{\alpha_1},\sigma_{\alpha_2},...,\sigma_{\alpha_n})$ and $(\sigma_{S_1},\sigma_{S_2},...,\sigma_{S_n})$ are in the same orbit of the action of $B_{n}$ .
\end{thm} 

To prove the theorem, we need the following definition.

\begin{defn}
	A sequence $E=(E_1,E_2,...,E_n)$ in $mod\Lambda$ is called a projective sequence if we have the following properties:
	
	For $1\leq r\leq n$, let $S(E(r))$ be the set of composition factors of $E_1,E_2,...,E_r$, and $\cc(S(E(r)))$ be the thick closure of all the simple objects appearing in $S(E(r))$. When $S(E(r))$ consists of simple objects, its thick closure is just its extension closure. 
	
	1)The number of simple objects appearing in $S(E(r))$ is $r$.
	
	2)$E_r$ is a projective object in $\cc(S(E(r)))$.
	
	3)$\top(E_r)\notin S(E(r))$.
	
\end{defn} 
Now we have the following lemma.

\begin{lem}
	A projective sequence is an exceptional sequence.
\end{lem} 

\begin{pf}
	Let $E=(E_1,E_2,...,E_n)$ be a projective sequence. Take a $1\leq r\leq n$. Since $E_r$ is a projective module in $\cc(E^r)$, $\Ext_\Lambda(E_r,E_j)=0$ for $1\leq j\leq r$. By property (3), $\top(E_r)\notin S(E(r-1))$, so $\Hom_\Lambda(E_r,E_j)=0$ for $1\leq j<r$. Then it follows that $E$ is an exceptional sequence.
	\end{pf}

\begin{lem}
	If $\omega\in W$ with $\lvert\omega\rvert_a=1$, then for each decomposition of $\omega=\sigma_{S_{i_1}}\sigma_{S_{i_2}}...\sigma_{S_{i_t}}$, we can delete some $\sigma_{S_{i_j}}$ such that $\sigma_{S_{i_1}}\sigma_{S_{i_2}}...\sigma_{S_{i_{j-1}}}\sigma_{S_{i_{j+1}}}...\sigma_{S_{i_t}}=1$, i.e. $\sigma_{S_{i_1}}\sigma_{S_{i_2}}...\hat{\sigma}_{S_{i_j}}...\sigma_{S_{i_t}}=1$
\end{lem} 

\begin{pf} Since $\lvert\omega\rvert_a=1$, we have $\omega=\sigma_{S_{i_1}}\sigma_{S_{i_2}}...\sigma_{S_{i_t}}=\sigma_{\alpha}$ for some real root $\alpha$. Since $\sigma_{\alpha}(\alpha)<0$, there must be some $\sigma_{S_{i_j}}$ such that $S_{i_j}=\sigma_{S_{i_{j+1}}}\sigma_{S_{i_{j+2}}}...\sigma_{S_{i_t}}(\alpha)$, the following holds:
$$\sigma_{\alpha}=\sigma_{S_{i_1}}\sigma_{S_{i_2}}...\sigma_{S_{i_{j-1}}}\sigma_{S_{i_j}}\sigma_{S_{i_{j+1}}}...\sigma_{S_{i_t}}$$$$=\sigma_{S_{i_1}}\sigma_{S_{i_2}}...\sigma_{S_{i_{j-1}}}\sigma_{S_{i_{j+1}}}...\sigma_{S_{i_t}}\sigma_{\alpha}\sigma_{S_{i_t}}...\sigma_{S_{i_{j+1}}}\sigma_{S_{i_{j+1}}}...\sigma_{S_{i_t}}$$$$=\sigma_{S_{i_1}}\sigma_{S_{i_2}}...\sigma_{S_{i_{j-1}}}\sigma_{S_{i_{j+1}}}...\sigma_{S_{i_t}}\sigma_{\alpha}$$
$$\sigma_{S_{i_1}}\sigma_{S_{i_2}}...\hat{\sigma}_{S_{i_j}}...\sigma_{S_{i_t}}$$
Thus the lemma is proved.  
\end{pf}

\begin{defn}
   For a reflection $\omega\in W$ at some real root, define $l'(\omega)$ to be the minimal length of $\omega_1$ such that $\omega=\omega_1\sigma_{S_j}\omega_1^{-1}$ for some simple reflection $\sigma_{S_j}$.
\end{defn}
The following lemma is well-known. For example, see \samepage\cite{Ringel}

\begin{lem}
	Given an exceptional sequence $S=(S_1,...,S_n)$ consisting of simple objects. The Coxeter element $c=\sigma_{S_1}\sigma_{S_2}...\sigma_{S_n}$ has the absolute length $n$.
\end{lem}
\begin{pf}
   Set $c=\sigma_{S_1}\sigma_{S_2}...\sigma_{S_n}=\sigma_{\alpha_1}\sigma_{\alpha_2}...\sigma_{\alpha_m}$ where $m=\lvert c\rvert_a$.	The lemma is a corollary of the following Lemma 4.12.
\end{pf}
\begin{lem}
   Given an exceptional sequence $S=(S_1,...,S_n)$ consisting of simple objects. If $\sigma_{S_{i_1}}\sigma_{S_{i_2}}...\sigma_{S_{i_t}}=\sigma_{\alpha_1}\sigma_{\alpha_2}...\sigma_{\alpha_s}$ ($i_1<i_2<...<i_t$) where $s=\lvert\sigma_{S_{i_1}}\sigma_{S_{i_2}}...\sigma_{S_{i_t}}\rvert_a$, then the following holds:
   
   a) There exists $(\sigma_{\beta_1},...,\sigma_{\beta_t})$ in the orbit of $(\sigma_{\alpha_1},...,\sigma_{\alpha_t})$ and $l'(\sigma_{\beta_i})$ is minimal in the $B_{s-i+1}$ action orbit of $(\sigma_{\beta_i},...,\sigma_{\beta_s})$. 
   
   b) there exists an $r$ such that $\sigma_{S_{i_1}}...\sigma_{S_{i_{r-1}}}\sigma_{S_{i_{r+1}}}...\sigma_{S_{i_t}}=\sigma_{\beta_1}\sigma_{\beta_2}...\sigma_{\beta_{s-1}}$. 
\end{lem}
\begin{pf}
	a) Let $\sigma_{\beta_1}$ be the element in the orbit of $(\sigma_{\alpha_1},...,\sigma_{\alpha_t})$ of the action of $B_n$ which has the minimal $l'(\cdot)$. Let $(\sigma_{\beta_1},\sigma_{\alpha_2^{(1)}},...,\sigma_{\alpha_t^{(1)}})$ be the element. Let $\sigma_{\beta_2}$ be the element in the orbit of $(\sigma_{\alpha_2^{(1)}},...,\sigma_{\alpha_t^{(1)}})$ of $B_{t-1}$ action that has the minimal $l'(\cdot)$. We do this process inductively. Finally we will get a sequence $(\sigma_{\beta_1},...,\sigma_{\beta_t})$ which  satisfies the required properties.
	
   b) For each $\sigma_{\beta_i}$, by the definition of $l'(\sigma_{\beta_i})$, there exists a $\omega_i$ such that $l(\omega_i)=l'(\sigma_{\beta_i})$ and $$\sigma_{\beta_i}=\omega_i\sigma_{S_{j_i}}\omega_i^{-1}$$ Then by the assumption in the lemma, we have
   $$\sigma_{\beta_s}=\sigma_{\beta_{s-1}}...\sigma_{\beta_1}\sigma_{S_{i_1}}...\sigma_{S_{i_t}}$$
   Now we apply Lemma 4.9. It is divided into four cases:
   
   Case1: If we delete some simple reflection $\sigma_{S_{j_i}}$, then by Lemma 4.9
   $$1=\sigma_{\beta_{s-1}}...\sigma_{\beta_{i+1}}\sigma_{\beta_{i-1}}...\sigma_{\beta_1}\sigma_{S_{i_1}}...\sigma_{S_{i_t}}$$
   which contradicts the assumption that $s=\lvert\sigma_{S_{i_1}}\sigma_{S_{i_2}}...\sigma_{S_{i_t}}\rvert_a $.
   
   Case2: If we delete some simple reflection in $\omega_i$. We write $\omega_i=\sigma_{S_{k_1}}...\sigma_{S_{k_p}}$ where $p=l(\omega)$. Suppose $\sigma_{S_{k_q}}$ is deleted. By Lemma 4.9
   $$1=\sigma_{\beta_{s-1}}...\sigma_{\beta_{i+1}}\sigma_{S_{k_1}}...\hat{\sigma}_{S_{k_q}}...\sigma_{S_{k_p}}\sigma_{S_{j_i}}\omega_i^{-1}\sigma_{F_{i-1}}...\sigma_{F_1}\sigma_{S_{i_1}}...\sigma_{S_{i_t}}
   $$
   $$=\sigma_{\beta_{s-1}}...\sigma_{\beta_{i+1}}\sigma_{S_{k_1}}...\sigma_{S_{k_{q-1}}}\sigma_{S_{k_{q+1}}}...\sigma_{S_{k_p}}\sigma_{S_{j_i}}\omega_i^{-1}\sigma_{F_{i-1}}...\sigma_{F_1}\sigma_{S_{i_1}}...\sigma_{S_{i_t}}$$
   $$=\sigma_{\beta_{s-1}}...\sigma_{\beta_{i+1}}\sigma_{S_{k_1}}...\sigma_{S_{k_{q-1}}}\sigma_{S_{k_q}}\sigma_{S_{k_{q-1}}}...\sigma_{S_{k_1}}\sigma_{S_{k_1}}...\sigma_{S_{k_q}}\sigma_{S_{k_{q+1}}}...\sigma_{S_{k_p}}\sigma_{S_{j_i}}\omega_i^{-1}\sigma_{\beta_{i-1}}...\sigma_{\beta_1}\sigma_{S_{i_1}}...\sigma_{S_{i_t}}$$
   $$=\sigma_{\beta_{s-1}}...\sigma_{\beta_{i+1}}\sigma_{\beta'}\sigma_{\beta_i}\sigma_{\beta_{i-1}}...\sigma_{\beta_1}\sigma_{S_{i_1}}...\sigma_{S_{i_t}}$$
   where $\sigma_{\beta'}=\sigma_{S_{k_1}}...\sigma_{S_{k_{q-1}}}\sigma_{S_{k_q}}\sigma_{S_{k_{q-1}}}...\sigma_{S_{k_1}}$
   
   But by $B_s$ action $(\sigma_{\beta_1},...,\sigma_{\beta_i},\sigma_{\beta'},...,\sigma_{\beta_s})$ and $(\sigma_{\beta_i},...,\sigma_{\beta_s})$ are in the same orbit. And by the expression of $\sigma_{\beta'}$ $l'(\sigma_{\beta'})<l'(\sigma_{\beta_i})$. Which is a contradiction.
   
   Case3: If we delete some simple reflection in $\omega_i^{-1}$. This is similar to Case2.
   
   Case4: If we delete some $\sigma_{S_{i_r}}$, it is exactly what the lemma says.

\end{pf}

Apply Lemma 4.12 b) inductively we can get $s=t$. So apply this to Lemma 4.11 then it holds.

Now we come to prove Theorem 4.6. 
\begin{pf}
  By Lemma 4.12 a), there exists $(\sigma_{\beta_1},...,\sigma_{\beta_n})$ in the orbit of $(\sigma_{\alpha_1},...,\sigma_{\alpha_n})$ and $l'(\sigma_{\beta_i})$ is minimal in the $B_{s-i+1}$ action orbit of $(\sigma_{\beta_i},...,\sigma_{\beta_n})$. Since $\sigma_{\beta_i}$ $1\leq i\leq n$ are real roots, there exists a unique indecomposable module $F_i$ satisfying $\udim F_i=\beta_i$.
  
   Now we have
   $$\sigma_{F_1}\sigma_{F_2}...\sigma_{F_n}=\sigma_{S_1}...\sigma_{S_n}$$
   
   By Lemma 4.8, what we need to prove firstly is that $F=(F_1,F_2,...,F_n)$ is a projective sequence.
   
   Denote $S(F(r))$ be the set of composition factors of $F_1,...,F_r$.
   
   By applying Lemma 4.12 inductively, the following properties of $F=(F_1,F_2,...,F_n)$ follow:
   
   1)The number of simple objects appearing in $S(F(r))$ is $r$.
   
   2)$\top(F_r)\notin S(F(r))$.
   
   By our assumption, $S=(S_1,...,S_n)$ is a complete exceptional sequence. If $S(F(r))=\{S_{i_1^{(r)}},...,S_{i_r^{(r)}}\}$ where $i_1^{(r)}<i_2^{(r)}<...<i_r^{(r)}$, $(S_{i_1^{(r)}},...,S_{i_r^{(r)}})$ is a complete exceptional sequence of $\cc(S(F(r)))$. By Lemma 4.12, $\udim(F_r)=\sigma_{S_{i_r^{(r)}}}\sigma_{S_{i_{r-1}^{(r)}}}...\sigma_{S_{i_k^{(r)}}}(\udim S_{i_{k-1}^{(r)}})$ for some $k$, $1\leq k\leq r$. Thus $F_r$ is a projective object in $\cc(S(F(r)))$. So $F$ is a projective sequence. By Lemma 4.8, $F$ is an exceptional sequence. By Theorem 3.9, $F$ and $S$ are in the same orbit. If $g\cdot S=F$ for some $g\in B_n$, by $\udim R_YX=\pm\sigma_Y(X)$,  
  $ \udim L_XY=\pm\sigma_X(Y)$, we get  $g\cdot(\sigma_{S_1},\sigma_{S_2},...,\sigma_{S_n})=(\sigma_{F_1},\sigma_{F_2},...,\sigma_{F_n})$. The theorem is proved.

\end{pf}

\section{Ingalls-Thomas-Igusa-Schiffler bijection}

In the previous section, we label the simple modules in an appropriate order such that $S$ is a complete exceptional sequence. In this section, we write $c(\Lambda)=\sigma_{S_1}...\sigma_{S_n}$. By then we have two posets: $\bA(mod \Lambda)$ and $Nc(W(\Lambda),c(\Lambda))$. 

By Theorem 4.6 and Theorem 3.9 we have

\begin{lem}
	If $S=(S_1,S_2,...,S_n)$ is the complete exceptional sequence consisting of simple modules, and for n exceptional modules $\{E_1,.E_2,...,E_n\}$ we have
	\[\sigma_{S_1}\sigma_{S_2}...\sigma_{S_n}=\sigma_{E_1}\sigma_{E_2}...\sigma_{E_n}\]
	then $E=(E_1,E_2,...,E_n)$ is a complete exceptional sequence.
\end{lem}

\begin{pf} By Theorem 4.6, $(\sigma_{S_1},...,\sigma_{S_n})$ and $(\sigma_{E_1},...,\sigma_{E_n})$ are in the same orbit under braid group action. So $$g\cdot(\sigma_{S_1},...,\sigma_{S_n})=(\sigma_{E_1},...,\sigma_{E_n})$$ 
for some $g\in B_n$. By the formulas in Lemma 3.10:
$$\udim R_YX=\pm\sigma_Y(X)$$
$$\udim L_XY=\pm\sigma_X(Y)$$
If $\rho_i(\sigma_{E'_1},...,\sigma_{E'_n})=(\sigma_{F'_1}...\sigma_{F'_n})$,  $E_i,F_j$ in $T^n$, then for modules $\rho_i(E'_1,...,E'_n)=(F'_1,...,F'_n)$ in the level of exceptional sequence. So $(E_1,...,E_n)=g^{-1}\cdot(S_1,...,S_n)$ is a complete exceptional sequence.
\end{pf}

Now we can prove as same as in \cite{Ringel}.

\textbf{Main Theorem (Ingalls-Thomas-Igusa-Schiffler)} The map
\[cox:\ \bA(mod\Lambda)\rightarrow Nc(W(\Lambda),c(\Lambda))\]
is a poest isomorphism.

\begin{pf} If $cox(\mathcal{A})$=$cox(\mathcal{B})$, then we can choose a complete exceptional sequence $A=(A_1,...,A_r)$($B=(B_1,...,B_t)$) for $\mathcal{A}$($\mathcal{B}$), then we must have
\[\sigma_{A_1}\sigma_{A_2}...\sigma_{A_r}=\sigma_{B_1}\sigma_{B_2}...\sigma_{B_t}\]

By Lemma 3.4, we can extend $A$ to $A'=(A_1,...,A_r,A_{r+1},...,A_n)$, so we can see that $\sigma_{B_1}\sigma_{B_2}...\sigma_{B_t}\sigma_{A_{r+1}}...\sigma_{A_n}$ is the Coxeter element. By Lemma 5.1, $(B_1,...,B_t,A_{r+1},...,A_n)$ is a complete exceptional sequence. So we have 
\[\mathcal{B}=(^{\perp}\mathcal{A})^{\perp}=\mathcal{A}\]
which proves the injection.

For surjection, if $\omega=\sigma_{E_1}\sigma_{E_2}...\sigma_{E_r}\in Nc(W(\Lambda),c(\Lambda))$, then we can extend it to $\sigma_{E_1}\sigma_{E_2}...\sigma_{E_n}$ which is the Coxeter element. By Lemma 5.1, $(E_1,E_2,...,E_n)$ is a complete sequence. We just need to choose $\mathcal{A}$=$\cc(E_1,...,E_r)$. 

Finally we prove that this is a poset morphism, i.e. $\ca\leq\cb\Leftrightarrow cox(\ca)\leq cox(\cb)$. If $\ca\leq\cb$, choose a complete exceptional sequence  $A=(A_1,...,A_r)$ of $\ca$. By Lemma 3.4, it can be extended to a complete exceptional sequence  $B=(A_1,...,A_r,B_{r+1},...,B_t)$ of $\cb$, thus $cox(\ca)\leq cox(\cb)$ by the definition. Conversely, if $cox(\ca)\leq cox(\cb)$, i.e. $\sigma_{A_1}\sigma_{A_2}...\sigma_{A_r}\leq\sigma_{B_1}\sigma_{B_2}...\sigma_{B_t}$ where $A=(A_1,...,A_r)$ ($resp.B=(B_1,...,B_t)$) is a complete exceptional sequence of $\ca$ ($resp.\cb$). So there exists $C_{r+1},...,C_t\in\cb$  such that $\sigma_{A_1}\sigma_{A_2}...\sigma_{A_r}\sigma_{C_{r+1}}...\sigma_{C_t}=\sigma_{B_1}\sigma_{B_2}...\sigma_{B_t}$. Thus $\ca\leq\cb$. We finish the proof.
\end{pf} 
\begin{rem}
   In fact, all the consequence above can be generalized to the Artin hereditary algebras. Given an Artin hereditary algebra $\Lambda$, we can define the Cartan matrix associated to $K_0(\Lambda-mod)$. It is a symmetrizable generalized Cartan matrix. So there is a Kac-Moody Lie algebra $\fkg(\Lambda)$. Like what we have done in this note, there is a poset isomorphism 
   $$cox:\bA(mod \Lambda) \longrightarrow Nc(W(\Lambda),c(\Lambda))$$
\end{rem}
\section*{Acknowledgments}
Prof. Jie Xiao first introduced me this topic. Prof. Bangming Deng and Prof. Jie Xiao suggested me Ringel's paper and referred the material that I would need. They have carried on the thorough discussion on this theorem with me. I really thank them for their help and guidance.

\bibliographystyle{plain}
\bibliography{moban}

\end{document}